\documentclass{amsart}
\usepackage{amsmath}
\usepackage{amsfonts}
\usepackage{graphics}
\usepackage{epsfig}
\usepackage{amssymb}
\usepackage{amscd}

\newtheorem{thm}{Theorem}

\newtheorem{prop}[thm]{Proposition}
\theoremstyle{definition}

\newtheorem*{ack}{Acknowledgement}
\theoremstyle{remark}

\theoremstyle{definition}

\def\rchi{{\hbox{\raise1.5pt\hbox{$\chi$}}}}
\def\Aut{{\text{\rm{Aut}}}}

\begin{document}

\title[Periods of Strebel Differentials]
{Periods of Strebel Differentials and
Algebraic Curves Defined over the Field of
Algebraic Numbers}
\author[M.~Mulase]{Motohico Mulase$^1$}
\address{
Department of Mathematics\\
University of California\\
Davis, CA 95616--8633}
\email{mulase@math.ucdavis.edu}
\author[M.~Penkava]{Michael Penkava$^{2}$}
\address{
Department of Mathematics\\
University of Wisconsin\\
Eau Claire, WI 54702--4004}
\email{penkavmr@uwec.edu}
\date{July 16, 2001}
\thanks{$^{1}$
Research supported by NSF Grant
DMS-9971371 and the University of California, Davis.}
\thanks{$^{2}$
Research supported in part by NSF Grant
DMS-9971371
and
  the University of Wisconsin-Eau Claire.}

\allowdisplaybreaks
\setcounter{section}{0}
\begin{abstract}
In \cite{MP} we have shown that if a compact Riemann surface
admits a Strebel differential with rational periods, then the
Riemann surface is the complex model of an algebraic curve defined
over the field of algebraic numbers. We will show in this article
that even if all geometric data are defined over
$\overline{\mathbb{Q}}$, the Strebel differential can still have a
transcendental period. We construct a Strebel differential $q$ on
an arbitrary complete nonsingular algebraic curve $C$ defined over
$\overline{\mathbb{Q}}$ such that (i) all poles of $q$ are
$\overline{\mathbb{Q}}$-rational points of $C$; (ii) the residue
of $\sqrt{q}$ at each pole is a positive integer; and
     (iii) $q$ has a transcendental period.
\end{abstract}
\maketitle

\section{Introduction}\label{sect: intro}
The periodic function $e^{i\pi \theta}$ has the following
remarkable property: for every rational number $\theta$, $e^{i\pi
\theta}$ gives an algebraic number, while for every non-rational
algebraic number $\theta$, it gives a transcendental number. It
follows that if $e^{i\pi \theta}$ is algebraic and $\theta$ is not
rational, then $\theta$ is transcendental. This somewhat
reciprocal algebraicity-transcendence relation between the periods
of a function and its values is commonly seen
     among the periodic functions appearing
in algebraic geometry such as modular functions
and theta functions \cite{Baker}.

The purpose of this paper is to examine a
     similar  relation
between the periods of Strebel differentials on
a complete nonsingular algebraic curve and
its field of definition. In an earlier paper
\cite{MP}, we have established that if a
compact Riemann surface $C$ admits a Strebel
differential $q$ such that the length of
every \emph{critical
horizontal trajectory} is rational,
     then $C$ is the complex model of an
algebraic curve defined over the field
$\overline{\mathbb{Q}}$ of algebraic numbers.
(We refer to  \cite{MP}
for the  definitions
of the notions used in this article.)
Let us call the length of a critical horizontal
trajectory a \emph{period} of $q$. Thus the
rationality of the periods makes the field of
definition of $C$ algebraic. Following the
analogy of the relation between the periods
and the values of theta functions, one can ask
questions about the transcendence of the periods
of a Strebel differential and the field of
definition of an algebraic curve.

Let $\Gamma$ be a \emph{ribbon graph}, that is, a
graph with
a prescribed cyclic order of half-edges at each
vertex. A ribbon graph is a \emph{metric} ribbon
graph if a positive real number, conventionally
referred to as the \emph{length}, is assigned to
each edge. From every metric ribbon graph
$\Gamma$, one can construct  a compact Riemann surface $C$ and a
\emph{Strebel differential}
     $q$ on it in a
unique manner. The number of
poles of $q$ is equal to
the number of \emph{boundary
circuits} of $\Gamma$, and the residue of
$\sqrt{q}$ at each pole is the total length of
the boundary circuit corresponding to the pole.

Conversely, starting with a compact Riemann
surface $C$ with $n>0$ marked points and an
$n$-tuple of positive real numbers, one constructs
a unique metric ribbon graph as the union of
the critical
horizontal trajectories of the Strebel differential
$q$ on $C$. The integral of $\sqrt{q}$ between
zeros of $q$ along a critical horizontal
trajectory is called a \emph{period} of $q$.
The periods
of the Strebel differential determine a metric on the ribbon graph.

Harer \cite{Harer} used this idea to
establish an orbifold isomorphism
\begin{equation}
\label{eq: graph=moduli}
\coprod_{\Gamma} \frac{\mathbb{R}_+ ^{e(\Gamma)}}
{\Aut_{\partial}(\Gamma)}
\overset\sim{\longrightarrow}
     \mathfrak{M}_{g,n} \times \mathbb{R}_+ ^n,
\end{equation}
where $\Gamma$ runs over all ribbon graphs with valence
of each vertex no less than 3 and with
Euler characteristic $2-2g-n$, $e(\Gamma)$
is the number of edges of $\Gamma$,
$\Aut_{\partial}(\Gamma)$ is the group of
ribbon graph automorphisms
     of $\Gamma$ that fix its
boundary circuits, and $\mathfrak{M}_{g,n}$
is the moduli space of smooth compact
complex algebraic curves
of genus $g$ with $n$ ordered marked points.
     This isomorphism plays a key role in many recent papers
(cf. \cite{BIZ},
\cite{Harer-Zagier}, \cite{Kontsevich},
\cite{Mulase1995},
\cite{Penner}, \cite{Witten}).
For an explicit construction of this isomorphism,
we refer to \cite{MP}.

The set of \emph{rational points}
\begin{equation*}
\coprod_{\Gamma} \frac{\mathbb{Q}_+ ^{e(\Gamma)}}
{\Aut_{\partial}(\Gamma)}
\subset
\coprod_{\Gamma} \frac{\mathbb{R}_+ ^{e(\Gamma)}}
{\Aut_{\partial}(\Gamma)}
\end{equation*}
is well defined (although
not as an orbifold over $\mathbb{Q}$) because a
ribbon graph automorphism acts as a permutation of
edges. Let
$\mathfrak{M}_{g,n}(\overline{\mathbb{Q}})$
denote the moduli space of $n$-pointed complete
nonsingular algebraic curves defined over
$\overline{\mathbb{Q}}$.
Comparing Belyi's theorem \cite{Belyi},
Grothendieck's idea of \emph{dessins d'enfants}
\cite{Schneps}, and  Strebel
theory \cite{Strebel},
we have established in \cite{MP}  that there is a
natural injective map $j$:
\begin{equation}
\label{eq: CD}
\begin{CD}
\coprod_{\Gamma} \frac{\mathbb{R}_+ ^{e(\Gamma)}}
{\Aut_{\partial}(\Gamma)}
@>{\sim}>>
     \mathfrak{M}_{g,n} \times \mathbb{R}_+ ^n\\
\bigcup @. \bigcup\\
\coprod_{\Gamma} \frac{\mathbb{Q}_+ ^{e(\Gamma)}}
{\Aut_{\partial}(\Gamma)}
@>>{j}> \mathfrak{M}_{g,n}(\overline{\mathbb{Q}})
\times \mathbb{Q}_+ ^n .\\
\end{CD}
\end{equation}
Belyi's theorem shows that \emph{every}
complete nonsingular algebraic curve over
$\overline{\mathbb{Q}}$ is constructed
as the image of the map $j$ \emph{if} we do not
specify the number $n$ of marked points.
However, in the light of the geometric
$Gal(\overline{\mathbb{Q}}/\mathbb{Q})$
actions \cite{Schneps}, $j$ cannot
be surjective. Indeed, we shall prove the following:
\begin{thm}
\label{th: main}
Let $C$ be an arbitrary complete nonsingular
algebraic curve defined over
$\overline{\mathbb{Q}}$. Then there is a
Strebel differential $q$ on the complex model
of $C$ such that
\begin{enumerate}
\item every pole of $q$ is a
$\overline{\mathbb{Q}}$-rational point of $C$;
\item the residue of $\sqrt{q}$ at each pole
is a positive integer; and
\item $q$ has a transcendental period.
\end{enumerate}
\end{thm}

\begin{figure}[htb]
\centerline{\epsfig{file=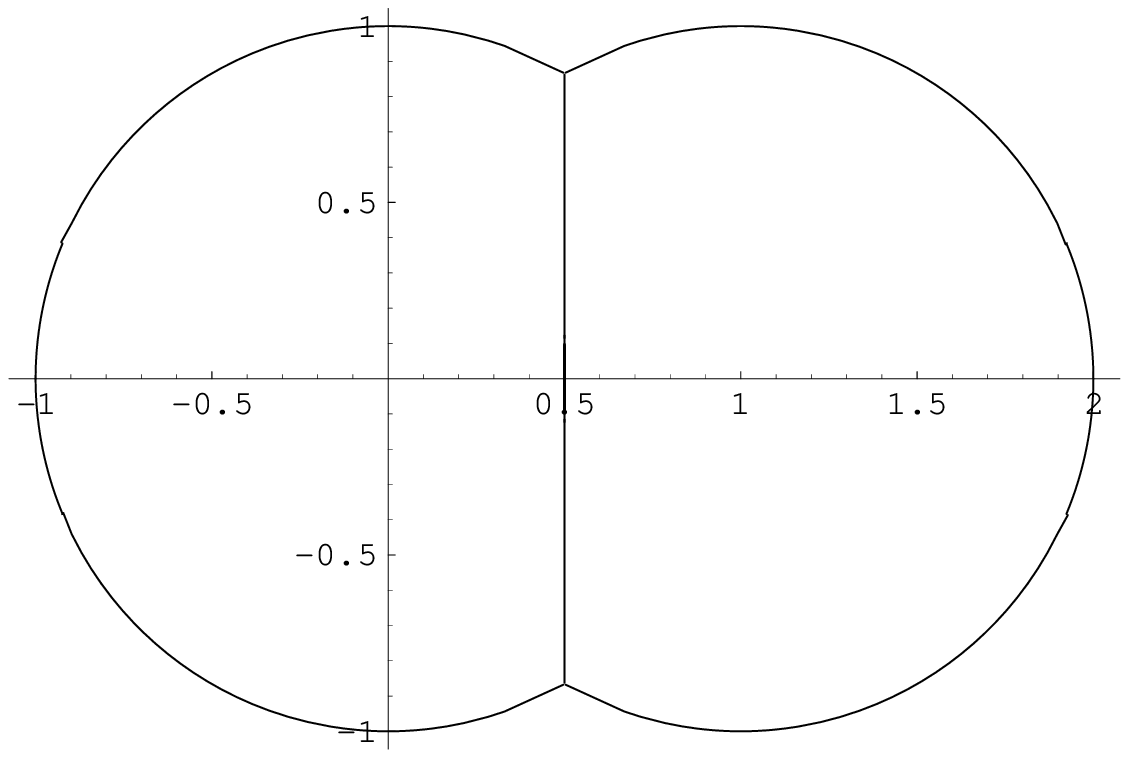, width=2in}}
\caption{Graph $\Gamma_y$, the
     inverse image of $[0,1]$
via $\zeta = \frac{4(y^2-y+1)^3}{27y^2(1-y)^2}$.}
\label{fig: Gammay}
\end{figure}

\begin{figure}[htb]
\centerline{\epsfig{file=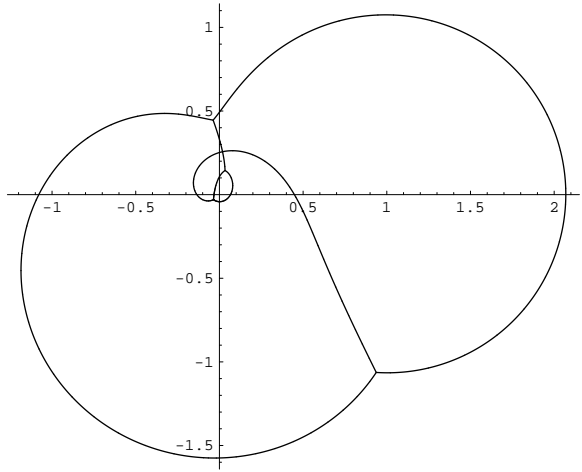, width=3in}}
\caption{Graph $\Gamma_c$ for
$s = \frac{5\sqrt{3}}{11}$.}
\label{fig: Gammac}
\end{figure}

\begin{ack}
The authors would like to thank Josephine Yu for
creating the graphics Figure~\ref{fig: Gammay} and
Figure~\ref{fig: Gammac} in this paper.
\end{ack}

\section{Construction of the Strebel Differential}
\label{sect: Strebel}

Let us start by constructing a simple example
on $\mathbb{P}^1$. We wish to exhibit a
Strebel differential $q_c$ on $\mathbb{P}^1$
that has poles at $0$, $1$, $\infty$, $c$ and
$c^2$, such that $\sqrt{q_c}$ has residues
    2, 2, 2, 4, and 2, resp., at these poles,
with $c\in \overline{\mathbb{Q}}\setminus \{0,1\}$
    a constant to be determined later.
We define two rational maps
$f$ and $g$ in order to construct certain Strebel
differentials on $\mathbb{P}^1$. First we choose
\begin{equation}
\label{eq: f}
f:\mathbb{P}^1\owns x
\longmapsto
y=\frac{1}{(1-c)^2}\; \frac{(x-c)^2}{x}\in
\mathbb{P}^1.
\end{equation}
Since
$$
f'(x)=\frac{1}{(1-c)^2}\;\frac{x^2-c^2}{x^2},
$$
$f$
is ramified at $x=\pm c$. We note that
$$
f(0)=f(\infty)=\infty,\quad f(1) = f(c^2)=1,
\quad f(c)=0, \quad f(-c)=-\frac{4c}{(1-c)^2}.
$$
The other rational map is
\begin{equation}
\label{eq: g}
g:\mathbb{P}^1 \owns y\longmapsto
\zeta = \frac{4(y^2-y+1)^3}{27y^2(1-y)^2}\in
\mathbb{P}^1.
\end{equation}
Consider the meromorphic quadratic differential
$$
q_0 = \frac{1}{4\pi^2}\;
\frac{d\zeta^2}{\zeta(1-\zeta)}
$$ on
$\mathbb{P}^1$, and put
\begin{equation}
\label{eq: q}
q_1=g^*(q_0)=-\frac{1}{\pi^2}\; \frac{y^2-y+1}{y^2(1-y)^2}dy^2.
\end{equation}
It has quadratic poles at $0$, $1$, and $\infty$,
and simple zeros at $\frac{1}{2} \pm i
\frac{\sqrt{3}}{2}$. The residue of $\sqrt{q_1}$
at each pole is $2$, and the three periods of $q_1$
are all $1$.
Let
$\Gamma_y=g^{-1}([0,1])$. It has been shown in
\cite{MP} that
$\Gamma_y$ is the ribbon graph consisting of the set
of critical trajectories of the Strebel differential
$q_1$.

Next, let
$s$ be a real algebraic number
such that $0<s<\frac{\sqrt{3}}{2}$,  choose the
constant
$c$ so that
\begin{equation}
\label{eq: c}
f(-c)=-\frac{4c}{(1-c)^2} = \frac{1}{2} + i s ,
\end{equation}
and define
$q_c =f^*(q_1)$.
Since the two critical values of the
double-sheeted holomorphic
covering map $f$ are
$f(c)=0$ and $f(-c)$,
one sees that the
Strebel differential
$q_c$ on
$\mathbb{P}^1$ has quadratic poles at
$0$, $1$, $\infty$, $c$ and $c^2$ with
residues of $\sqrt{q_c}$ at these poles
2, 2, 2, 4, and 2, resp., and a new double
zero at $-c$.
Let
$\Gamma_c = f^{-1}(\Gamma_y)$.  
Because of the choice of $c$ in Eq.~\ref{eq: c}, 
$\Gamma_c$ is the ribbon
graph consisting of the critical trajectories for the
Strebel differential $q_c$. It has  four tri-valent
vertices at $f^{-1}(\frac{1}{2} \pm i
\frac{\sqrt{3}}{2})$ and a unique 4-valent vertex at
$-c = f^{-1}(\frac{1}{2} + i s)$.

Let $L$ denote the edge of
$\Gamma_c$ that connects the vertex $-c$ and
one of the tri-valent vertices. From
Eq.~\ref{eq: q} one sees that the period of $q_c$
corresponding to the edge $L$ is equal to
$\frac{1}{2} \pm
\ell(s)$ depending on the choice of the
tri-valent vertex, where
\begin{equation}
\begin{aligned}
\label{eq: ell}
\ell(s) &= \int_{1/2} ^{1/2 + is}
\sqrt{-\frac{1}{\pi^2}\;
\frac{y^2-y+1}{y^2(1-y)^2}}dy\\
&= \frac{1}{\pi}\int_0 ^s
\frac{\sqrt{3/4 - t^2}}{1/4 + t^2}dt\\
&=\frac{1}{\pi}\;
\arcsin\left(
\frac{2s(9+4s^2)}{3\sqrt{3}(1+4s^2)}
\right).
\end{aligned}
\end{equation}
In Section~\ref{sect: trans}, we shall show that
$\ell(s)$ is transcendent at $s=r\sqrt{3}$
for every rational number $r$ such that
$0<r<\frac{1}{2}$.

The above construction immediately gives the
construction of the Strebel differential
$q$ for the general case of
Theorem~\ref{th: main}. Let $C$
be an arbitrary complete nonsingular
algebraic curve defined over
$\overline{\mathbb{Q}}$. Belyi \cite{Belyi}
has shown that there is a holomorphic map
\begin{equation}
\label{eq: beta}
\beta: C\longrightarrow \mathbb{P}^1
\end{equation}
that is ramified only over $0$, $1$ and
$\infty$. The map $\beta$
is called a \emph{Belyi map}.
Without loss of generality, we can assume
that the ramification degrees over the point
$0\in \mathbb{P}^1$ are
no less than 3 and the ramification degrees
over $1\in \mathbb{P}^1$ are always 2.
Define $q=\beta^* (q_c)$ and $\Gamma =
\beta^{-1}(\Gamma_c)$. Then $q$ is the Strebel
differential on $C$ with poles at
$\beta^{-1}(0)$, $\beta^{-1}(1)$,
$\beta^{-1}(\infty)$, $\beta^{-1}(c)$,
and $\beta^{-1}(c^2)$. The residue of
$\sqrt{q}$ at each of these poles is
2 times the ramification degree of $\beta$
at each ramification point, 4 at each
inverse image of $c$, and 2 at each inverse
image of $c^2$.
Certainly the period of $q$ corresponding to
any edge of $\Gamma$ that is an inverse image
of $L$ is a transcendental number if we chose
$s=r \sqrt{3}$.

\section{Transcendence of the Period}
\label{sect: trans}

Let us now show that $\ell(s)$ of
Eq.~\ref{eq: ell} is transcendental for
every $s=r\sqrt{3}$, where $r$ is a
rational number in between $0$ and $1/2$.
Since
$$
\sin(\pi \ell(s)) =
\frac{2s(9+4s^2)}{3\sqrt{3}(1+4s^2)}
=\frac{2r(3 +4r^2)}{(1+12r^2)}\in \mathbb{Q}
\setminus\{1/2\},
$$
the claim follows from
\begin{prop}
Let $0 < \ell <1/2$ be such that
$a = \sin(\pi\ell)$ is rational but not equal
to $\frac{1}{2}$.
Then $\ell$ is transcendental.
\end{prop}

\begin{proof}
Let
$$
b = -a + i\sqrt{1-a^2} = e^{i\pi(\ell +1/2)}.
$$
This is a solution of the quadratic equation $x^2+2ax +1=0$. Let
$a=m/n$ be an irreducible fraction. Then $nx^2+2mx +n$ is
primitive and irreducible  if $n$ is odd, and if $n=2k$, then
$kx^2 + mx + k$ is primitive and irreducible.

Suppose that $\ell\in\mathbb{Q}$. Then
$b$ is a solution of the primitive equation
$x^N - 1 = 0$ for some integer $N$. We note
that the primitive minimal polynomial of
$b$ divides $x^N-1$. Therefore, if $n$ is
odd, then $n=1$, and hence $a = m\ge 1$, which
is a contradiction. If $n$ is even, then
$k=1$ and the only possibility is $a=
\frac{1}{2}$. Thus $\ell$ is not rational.
It cannot be algebraic because $b$ is
algebraic. This completes the proof.
\end{proof}

It would be  desirable to establish
that if a compact Riemann surface $C$ admits
a Strebel differential $q$ whose periods are
algebraic but not rational, then the
geometric data $(C, (p_1, \cdots, p_n))$,
where $(p_1, \cdots, p_n)$ are the poles of
$q$, cannot be defined over
$\overline{\mathbb{Q}}$. In our simple example
of Section~\ref{sect: Strebel}, if the period
corresponding to the edge $L$ is taken to be
algebraic but not rational, then the ramification
points $c$ and $c^2$ become transcendental.
However, we do not have any general theorem
in this direction.

\providecommand{\bysame}{\leavevmode\hbox to3em{\hrulefill}\thinspace}

\bibliographystyle{amsplain}

\end{document}